# Sequences and Polynomial Congruence

Darrell Cox[a], Sourangshu Ghosh[b], Eldar Sultanow[c]

[a]Grayson County College
Denison, TX 75020
USA
dilycox@cableone.net

[b]Department of Civil Engineering,
Indian Institute of Technology Kharagpur, West Bengal, India
sourangshug123@gmail.com

[c]Potsdam University,
14482 Potsdam,
Germany
Eldar.Sultanow@wi.uni-potsdam.de

## Abstract

In this paper we shall find a new connection between $n^{th}$ degree polynomial mod p congruence with n roots and higher order Fibonacci and Lucas sequences. We shall first discuss about the recent work been done in sequences and their connection to polynomial congruence and then find out new relations between particular recurrence relation and the congruence of the sequences.

## 1. Introduction

Let the $U$ series is defined by the recurrence relation

$$U_i - a_1 U_{i-1} - a_2 U_{i-2} - \ldots - a_n U_{i-n} = 0$$

Here $i = 1, 2, 3, \ldots, U_0 = 1$, and $U_i = 0$ if $i < 0$.

In this article we shall consider the sequence modulo some prime number $p$. Now we shall start looking into some properties of polynomial congruence. Sun[1] used general linear recurring sequences by defining a recurrence relation to extend the Lucas Sequence. He established the following important condition for polynomials to be congruent to the linear polynomials product $(mod\ p)$ which we define as follows:

**Theorem 1.1[Sun[1]]:** Let $m \geq 2, m, a_1, \ldots, a_m \in Z$, and let $p$ be a prime such that $p > m$ and $p$ is coprime to $a_m$. Then the congruence $x_m + a_1 x_{m-1} + \cdots + a_m \equiv 0\ (mod\ p)$ has $m$ distinct solutions if and only if it satisfies the following condition:

$$U_{p-m} \equiv \cdots \equiv U_{p-2} \equiv 0\ (mod\ p) \quad and \quad U_{p-1} \equiv 1\ (mod\ p)$$

Sun[2] then used the Chebotarev density theorem which states that if the recurrence polynomial $f(x)$ is irreducible over $z(x)$ then the set $S$ of primes $p$ which will satisfy the condition of $f(x) \equiv 0\ (mod\ p)\ p \in S$ having $m$ solutions will have a positive density $d(S)$ or mathematically speaking the following holds:

$$d(S) = \lim_{x \to +\infty} \frac{|\{p:\ p \leq x, p \in S\}|}{|\{p:\ p \leq x, p\ is\ a\ prime\}|} > 0$$

To prove that if $f(x)$ over $Z(x)$ is irreducible, then there are must exist infinitely many prime $p$ satisfying the above stated condition.

One question of interest is when that sequence will repeat itself or what will be the period of that sequence. Hardy and Wright[2] found some upper limits of the period length of the Fibonacci by analyzing the properties of the roots of the recurrence polynomial. In this article we shall build on this approach to work in analyzing the recurrence relation roots and congruence of the sequences. Let us now define the term period for the sequence, we say that a sequence has period $l$ if $l$ is the smallest integer such that $U_{l+1} \equiv U_1 \pmod{p}$.

We shall now state and prove a theorem originally stated by Lehman[3] and Wallis[4] which puts a condition for the periodicity of recursive sequences modulo $p$.

**Theorem 1.2 [Lehman[3]]:** Let the $U$ series is defined as stated above by the recurrence relation $U_i - a_1 U_{i-1} - a_2 U_{i-2} - \ldots - a_n U_{i-n} = 0$. If $gcd(a_k, m) = 1$, then the sequence is periodic modulo $m$.

**Proof:** $Z_m$ will have exactly $m^k$ distinct $k-tuples$. Therefore there must exist two distinct integers $s$ and $t$ in $Z_m$ ( $0 \leq s < t \leq m^k$) such that it will satisfy $U_{s+i} \equiv U_{t+i}$ by the pigeonhole principle for $0 \leq i \leq k-1$. By the well-ordering principle there must exist a smallest non-negative integer $s > 0$ in $Z_m$.

If we take $s > 0$, then the condition $U_{s+k-1} \equiv U_{t+k-1}$ shall imply $a_k U_{s-1} \equiv a_k U_{t-1}$ since

$$a_1 U_{s+k-2} + a_2 U_{s+k-3} + \cdots + a_{k-1} U_s + a_k U_{s-1} = a_1 U_{t+k-2} + a_2 U_{t+k-3} + \cdots + a_{k-1} U_t + a_k U_{t-1}$$

Now If we have $gcd(a_k, m) = 1$ as the condition mentioned above in the theorem then we must have $a_{s-1} \equiv a_{t-1}$ in $Z_m$. This contradicts our initial assumption, therefore s must 0.

Therefore the sequence is periodic modulo $m$.

We shall now define the characteristic polynomial. If $U$ is a recursive sequence defined by the recurrence relation $U_i + a_1 U_{i-1} + a_2 U_{i-2} + \ldots + a_n U_{i-n} = 0$ as stated earlier then the characteristic polynomial of that sequence is defined to be

$$f(x) = x^i + a_1 x^{i-1} + a_2 x^{i-2} + \ldots + a_n x^{i-n}$$

The series can be written as the combination of the powers of the solutions of $f(x)=0$, whose values are determined by the sequence's initial terms.

Lehman[3] showed that the recursive sequence modulo some integer $m$ will be equal to the order of the formal root $\omega$ of the characteristic polynomial under some restrictions of $m$. Lehman[3] further showed using the properties of ring automorphisms of $Z_p[\omega]$ that if the recurrence polynomial $f$ has no repeated factors in the polynomial ring $Z_p[x]$ and $t$ is the LCM of the all the degrees of the irreducible factors of $f$ in $Z_p[x]$ then it is the smallest positive integer for which the order of the root $\omega$ shall divide $p^t - 1$. We shall state and prove the results now

Let us consider $\langle f \rangle$ to be the principal ideal of $Z_m[x]$ generated by the characteristic polynomial $f$. Consider the quotient ring $Z_m[x]/\langle f \rangle$. Here we define $Z_m[x]$ to be the ring of polynomials whose coefficients are in $Z_m$. Notice that every polynomial $g$ can be expressed in terms of 2 other unique polynomials $q$ and $r$ in $Z_m[x]$ such that it satisfies $g = f.q + r$, as r is the remainder it must be of smaller degree than $f$. Then we must have the two coset being $g + \langle f \rangle = r + \langle f \rangle$. Let us now define the ring $\frac{Z_m[x]}{\langle f \rangle}$ with the ring $Z_m[\omega]$ as

$$Z_m[\omega] = \{b_{k-1}\omega^{k-1} + b_{k-2}\omega^{k-2} + \cdots + b_1\omega^1 + b_0 \text{ and } \omega^k = \sum_{i=1}^{k} a_i \omega^{k-i}\}$$

Here we have written the coset $x + \langle f \rangle$ as $\omega$. We can now form a connection between recursive sequences with characteristic polynomial $f$ and the quotient ring $Z_m[x]/\langle f \rangle$. Let the $U$ series is defined by the recurrence relation $U_i - a_1 U_{i-1} - a_2 U_{i-2} - \ldots - a_n U_{i-n} = 0$ as defined earlier or $U(j,n) = \sum_{i=j}^{k} a_i U_{n-i}$. Then for all $n \geq k$ we have the following 2 conditions:

$$U(k,n) = a_k U_{n-k}$$

And $$U(j+1,n) + a_j U_{n-j} = U(j,n) \quad if\ 1 \leq j < k$$

**Theorem 1.3[Lehman³]:** Let $U$ be defined recursively as stated earlier. Then for every integer $n \geq 0$,

$$\alpha \omega^n = U_{n+k-1} \omega^{k-1} + \sum_{j=2}^{k} U(j, n+k+j-2) \omega^{k-j}$$

Here $\alpha$ is an element of $Z_m[\omega]$ which can be defined in terms of initial terms of the recursive sequence.

$$\alpha = U_{k-1} \omega^{k-1} + U(2,k) \omega^{k-2} + U(3,k) \omega^{k-3} + U(k-1, 2k-3)\omega + U(k, 2k-2)$$

From this theorem Lehman³ inferred that if $gcd(a_k, m) = 1$, and $Z_m[\omega] = Z_m[\omega]/\langle f \rangle$, then the sequence $U$ will be periodic modulo $m$ if and only if $\alpha \omega^l = \alpha$ in $Z_m[\omega]$. From this we can write a corollary that is the same as corollary 5 of Lehman³.

**Corollary 1.1[Lehman³]:** The sequence U is periodic modulo m, with the period $l$ dividing $ord_m(\omega)$, the order of $\omega$ in the group, $Z_m[\omega]$ with the operation $\times$, of units in $Z_m[\omega]$. If $A = \{\beta \in Z_m[\omega] | \alpha\beta = 0\}$, then $l$ is the order of $\omega + A$ in the group of units of the quotient ring $Z_m[\omega]/A$.

This recursive sequences are used to develop algorithms for factorization of higher order polynomials modulo primes by many authors (Adams[4], Lehman[3] and Sun[1]). There are some excellent works on recursive sequences modulo some positive integer. Some of them are Lehman[3], Engstrom[5], Ward[6] and Fillmore[7]. For the factorization of a polynomial into linear factors modulo a prime number $p$, Sun[4] and Skolem[8] find a condition which is basically the same as that of Corollary 13 of Lehman[3].

## 2. New relations between the recurrence relation and the sequences

Let $p$ be a prime greater than 2, $n$ a natural number such that $n \leq p - 2$, and $a_1, a_2, \ldots a_n$ integers. Let the U series be defined by the recurrence relation $U_i + a_1 U_{i-1} + a_2 U_{i-2} + \ldots + a_n U_{i-n}$ where $i = 1,2,3,\ldots, U_0 = 1$, and $U_i = 0$ if $i < 0$. The principal subject here is the proof of the following theorem;

**Theorem 2.1:** The congruence $x^n + a_1 x^{n-1} + a_2 x^{n-2} + \ldots + a_n \equiv 0 (mod\ p), 0 < x < p$, has $n$ roots if and only if $U_{p-1+i-n} \equiv U_{i-n} (mod\ p), i = 1,2,3,\ldots$.

(Also, if $a_1, a_2, \ldots a_n$ are elements of a Galois field, then the equation $x^n + a_1 x^{n-1} + a_2 x^{n-2} + \ldots + a_n \equiv 0 (mod\ p)$, has $n$ roots if and only if $U_{O-1+i-n} \equiv U_{i-n}$ where $O$ is the order of the Galois field.)

Only $n$ where $n \leq p - 2$ are considered since a mod p congruence of arbitrary degree can be reduced to degree $p - 2$ or less by using Fermat's theorem that $x^{p-1} \equiv 1 (mod\ p)$ if $p$ does not divide $x$. As will be shown, If $U_{p-1+i-n} \equiv U_{i-n}\ (mod\ p), i = 1,2,3, \ldots n$ then $U_{p-1+i-n} \equiv U_{i-n}\ (mod\ p), i = 1,2,3, \ldots$. This result combined with Theorem (2.1) gives a practical means of determining for which p a given nth degree congruence has n roots. (The same $U$ series applies for all p such that $n \leq p - 2$.) For example, some $U$ values corresponding to the congruence $x^3 - x^2 - x + 2 \equiv 0 (mod\ p), 0 < x < p$, are;

$$U_1 = 1 \quad U_2 = 2 \quad U_3 = 1 \quad U_4 = 1 \quad U_5 = -2 \quad U_6 = -3 \quad U_7 = -7 \quad U_8 = -6$$

$$U_9 = -7 \quad U_{10} = 1 \quad U_{11} = 6 \quad U_{12} = 21 \quad U_{13} = 25 \quad U_{14} = 34 \quad U_{15} = 17 \quad U_{16} = 1$$

$U_{16} \equiv 1 (mod\ 17)$ and $U_{15} \equiv U_{14} \equiv 0 (mod\ 17)$, Therefore this congruence has three roots if $p = 17$. Also, this congruence does not have three roots for any p such that $5 \leq p < 17$.

## 2.1 Relationship of U Series to Fibonacci Numbers

Notice that if $n = 2$ and $a_1 = a_2 = -1$, then $U_i$ is the $(i + 1)$th Fibonacci number. From the perspective here, the Fibonacci numbers are indexed improperly. Alternately, the $i$th Fibonacci number can be defined to be equal to $(\zeta_1^i - \zeta_2^i)/(\zeta_1 - \zeta_2)$ where $\zeta_1, \zeta_2$ are the roots of the equation $\zeta^2 - \zeta - 1 = 0$. A well known result is that p divides the $(p-1)$th Fibonacci number if and only if 5 (the discriminant of the equation $\zeta^2 - \zeta - 1 = 0$ is a quadratic residue of p.

The congruence equation $x^2 - x - 1 \equiv 0\ (mod\ p), 0 < x < p$, has two roots if and only if 5 is a quadratic residue of $p$. Theorem (2.1) is a generalization of these results. There also exists an alternate $U_i$ definition;

**Theorem 2.2:** $U_i$, $i \geq 0$, equals $\sum \zeta_1^{e_1} \zeta_2^{e_2} \ldots \zeta_n^{e_n}$ where $\zeta_1, \zeta_2, \ldots \zeta_n$ are the roots of $\zeta^n + a_1\zeta^{n-1} + a_2\zeta^{n-2} + \ldots + a_n \equiv 0$ and the summation is over all combinations of non-negative integers $e_1, e_2, \ldots e_n$ such that $e_1 + e_2 + \ldots + e_n = i$.

As will be shown, $a_1(U_{p-1} - 1) + 2a_2 U_{p-2} + 3a_3 U_{p-3} + \ldots + na_n U_{p-n} \equiv 0 (mod\ p)$. Therefore if $n > 1$ and $p$ does not divide $a_1, a_2, \ldots a_n$, the congruence to zero of any group of $n - 1$ of $U_{p-1} - 1, U_{p-2}, U_{p-3}, \ldots, U_{p-n}$ implies the congruence to zero of all of $U_{p-1} - 1, U_{p-2}, U_{p-3}, \ldots, U_{p-n}$. This result is of special significance in the case $n = 2$. If $n = 2$, a property of the U series is $U_{p-1} \equiv (a_1^2 - 4a_2)^{(p-1)/2} (mod\ p)$ (since if $n = 2$, $U_{p-1} = (\zeta_1^p - \zeta_2^p)/(\zeta_1 - \zeta_2)$, $(\zeta_1^p - \zeta_2^p)/(\zeta_1 - \zeta_2) \equiv (\zeta_1 - \zeta_2)^p/(\zeta_1 - \zeta_2)(mod\ p)$, and $(\zeta_1 - \zeta_2)^2 = [-(\zeta_1 + \zeta_2)]^2 - 4\zeta_1\zeta_2 = a_1^2 - 4a_2)$. Therefore $U_{p-1} \equiv 1 (mod\ p)$ where $n = 2$ if and only if $a_1^2 - 4a_2$ is a quadratic residue of p. If $U_{p-1} \equiv 1 (mod\ p)$ and $U_{p-2} \equiv 0 (mod\ p)$ where n=2, then p does not divide $a_2$ (since by the recurrence relation, $U_{p-2} = (-a_1)^{p-1} + k_1 a_2$ and $U_{p-1} = (-a_1)^{p-1} + k_2 a_2$ where $k_1, k_2$ are integers). Conversely, if $a_1^2 - 4a_2$ is a quadratic residue of $p$ and $p$ does not divide $a_2$ where $n = 2$, then $U_{p-1} \equiv 1 (mod\ p)$ and $U_{p-2} \equiv 0 (mod\ p)$ (since $a_1(U_{p-1} - 1) + 2a_2 U_{p-2} \equiv 0 (mod\ p)$). Therefore $U_{p-1} \equiv 1, U_{p-2} \equiv 0 (mod\ p)$ where $n = 2$ if and only if $a_1^2 - 4a_2$ is a quadratic residue of $p$ and $p$ does not divide a2. This shows the equivalence of Theorem (2.1) in the case $n = 2$ with the familiar result that the congruence $x^2 + a_1 x + a_2 \equiv 0 (mod\ p), 0 < x < p$, has two roots if and only if $a_1^2 - 4a_2$ is a quadratic residue of $p$ and $p$ does not divide $a_2$.

## 2.2 Relationship of U Series to Lucas' Series

Denote $\sum \zeta_1^{e_1} \zeta_2^{e_2} \ldots \zeta_n^{e_n}$ where $\zeta_1, \zeta_2, \ldots \zeta_n$ are the roots of $\zeta^n + a_1 \zeta^{n-1} + a_2 \zeta^{n-2} + \ldots + a_n \equiv 0$ and the summation is over all combinations of non-negative integers $e_1, e_2, \ldots e_n$ such that $e_1 + e_2 + \ldots + e_n = i$ .and exactly h of $e_1, e_2, \ldots e_n$ are non-zero by $V_{i,h}$. If $i < h$, let $V_{i,h} = 0$. The following theorem is also proved;

**Theorem 2.3:** The congruence $x^n + a_1 x^{n-1} + a_2 x^{n-2} + \ldots + a_n \equiv 0 (mod\ p), 0 < x < p, p$ does not divide $a_1 a_2 \ldots a_n$, has $n$ roots if and only if $V_{p,i} \equiv V_{1,i} (mod\ p), i = 2, 3, 4, \ldots, n$.

As will be shown, $V_{p,i}$ is always congruent $(mod\ p)$ to $V_{1,i}$ when $i = 1$. Lucas[10] denoted $(x_1^i - x_2^i)/(x_1 - x_2)$, $i = 1, 2, 3, \ldots$, where $x_1, x_2$ are the roots of $x^2 - Px + Q = 0$ and $P, Q$ are relatively prime integers by $u_i$ and $x_1^i + x_2^i$ by $v_i$. If n=2 and a1, a2 are relatively prime, then $U_0, U_1, U_2, \ldots$ is Lucas' $u_1, u_2, u_3, \ldots$ and $V_{1,1}, V_{2,1}, V_{3,1}$ is Lucas' $v_1, v_2, v_3$ The $U_i, V_{1,i}$ series can then be considered generalizations of Lucas' $u_i, v_i$ series. For typographical convenience, let $c(i, j)$ denote $i$ "choose" $j$ (a binomial coefficient). The $U$ and $V$ series are related as follows;

**Theorem 2.4:** $a_j U_i = (-1)^j \sum c(j+k,j) V_{i+j,j+k}, j = 1,2,3,\ldots,n$ where the summation is from $k = 0$ to $k = n - j$.

**Theorem 2.5:** $V_{i,j} = (-1)^j \sum c(j+k,j) a_{j+k} U_{i-j-k}, j = 1,2,3,\ldots,n$ where the summation is from $k = 0$ to $k = n - j$.

## 2.3  König's Theorem

The proof of Theorem (2.1) is based on König's theorem (proposed by Julius König[11] in 1882). König's theorem is this; Let $f(x) = c_0 x^{p-2} + c_1 x^{p-3} + c_2 x^{p-4} + \ldots + c_{p-2}$ where the $c's$ are integers and $c_{p-2}$ is not divisible by the prime p. Then $f(x) \equiv 0 \pmod{p}$ has real roots if and only if the determinant of the cyclic matrix

$$\begin{matrix}
c_0 & c_1 & c_2 & \ldots & c_{p-3} & c_{p-2} \\
c_1 & c_2 & c_3 & \ldots & c_{p-2} & c_0 \\
c_2 & c_3 & c_4 & \ldots & c_0 & c_1 \\
& & & \cdot & & \\
c_{p-2} & c_0 & c_1 & \ldots & c_{p-4} & c_{p-3}
\end{matrix}$$

is divisible by p. Denote this matrix by $C$. In order that it has at least k distinct real roots it is necessary that all $p - k$ rowed minors of $C$ be divisible by $p$. If also not all $p - k - 1$ rowed minors are divisible by p, the congruence has exactly k distinct real roots. Kronecker's[12] version of König's theorem will also be used in the proof. Kronecker's version is this; $f(x) \equiv 0 \pmod{p}$ has exactly $k$ roots if and only if the rank of $C$ is exactly $p - 1 - k$.

Fermat's theorem is essential to the formulation of König's theorem. If $f(x) \equiv 0 \pmod{p}$ has a root, then this root is also a root of $x^{p-1} \equiv 1 \pmod{p}$ and hence p divides the determinant of the resultant of $f(x)$ and $x^{p-1} - 1$, i.e., the cyclic matrix $C$. This is the motivation for part of the $U$ series definition; the $U$ series has been defined so that the "resultant" of $U_i + a_1 U_{i+n} + a_2 U_{i+n-1} + \ldots + a_n U_i$ and $U_{p-1+i} - U_i, i = 0,1,2,\ldots,(p-2)$, equals the resultant of $x^n + a_1 x^{n-1} + a_2 x^{n-2} + \ldots + a_n$ and $x^{p-1} - 1$. The condition $U_{p-1+i-n} \equiv U_{i-n} \pmod{p}$, $i = 1,2,3,\ldots$, is then the analogue of Fermat's theorem.

## 2.4  Proof of Theorem (2.1) (the Part Using Kronecker's Theorem)

First suppose that $U_{p-1+i-n} \equiv U_{i-n} \pmod{p}$, $i = 1,2,3,\ldots$, Denote the $(p-1)*(p-1)$ cyclic matrix

$$\begin{matrix}
a_2 & a_3 & a_4 & \ldots & a_n & 0 & \ldots & 0 & 0 & 1 & a_1 \\
a_3 & a_4 & a_5 & \ldots & 0 & 0 & \ldots & 0 & 1 & a_1 & a_2 \\
a_4 & a_5 & a_6 & \ldots & 0 & 0 & \ldots & 1 & a_1 & a_2 & a_3
\end{matrix}$$

$$\begin{pmatrix} \vdots & & & & \vdots & & \cdots & & \vdots & \vdots \\ 1 & a_1 & a_2 & \cdots & a_{n-2} & a_{n-1} & \cdots & 0 & 0 & 0 & 0 \\ a_1 & a_2 & a_3 & \cdots & a_{n-1} & a_n & \cdots & 0 & 0 & 0 & 1 \end{pmatrix}$$

by $A$. Then since $U_i + a_1 U_{i-1} + a_2 U_{i-2} + \ldots + a_n U_{i-n} = 0$ and $U_{p-1+i-n} \equiv U_{i-n} \pmod{p}$, $i = 1, 2, 3, \ldots (p-1)$ $A(U_{p-2}, U_{p-3}, U_{p-4}, \ldots, U_0) \equiv (0, 0, 0, \ldots, 0) \pmod{p}$. $U_{p-1+i-n} \equiv 0 \pmod{p}, i = 1, 2, 3, \ldots, (n-1)$, therefore $M_0(U_{p-n-1}, U_{p-n-2}, U_{p-n-3}, \ldots, U_0) \equiv (0, 0, 0, \ldots, 0) \pmod{p}$ where $M_0$ is the $(p-1) * (p-n)$ matrix obtained from $A$ by deleting its first $n-1$ columns. Denote the $(p-1) * (p-n-1)$ matrix by $R$.

$$\begin{pmatrix} 0 & 0 & 0 & \cdots & 0 & 0 & 0 & 1 \\ 0 & 0 & 0 & \cdots & 0 & 0 & 1 & a_1 \\ 0 & 0 & 0 & \cdots & 0 & 1 & a_1 & a_2 \\ 0 & 0 & 0 & \cdots & 1 & a_1 & a_2 & a_3 \\ & & & \vdots & & & & \\ a_{n-2} & a_{n-1} & a_n & \cdots & 0 & 0 & 0 & 0 \\ a_{n-1} & a_n & 0 & \cdots & 0 & 0 & 0 & 0 \\ a_n & 0 & 0 & \cdots & 0 & 0 & 0 & 0 \end{pmatrix}$$

$R(-U_{p-n-1}, -U_{p-n-2}, U_{p-n-3}, \ldots, -U_1) \equiv (a_1, a_2, a_3, \ldots, 0, 0, 1) \pmod{p}$ (since $U_{p-1} \equiv 1, U_{p-2} \equiv U_{p-3} \equiv U_{p-4} \equiv \ldots \equiv U_{p-n} \equiv 0 \pmod{p}$) and hence the last column of $M_0$ is linearly dependent $\pmod{p}$ on the other columns of $M_0$. Similarly, $M_0 (-U_{p-n}, -U_{p-n-1}, -U_{p-n-2}, \ldots, -U_1) \equiv (a_2, a_3, a_4, \ldots, 0, 1, a_1) \pmod{p}$ ( since $U_p \equiv U_1, U_{p-1} \equiv U_0, U_{p-2} \equiv U_{p-3} \equiv U_{p-4} \equiv \ldots \equiv Up - n + 1 \equiv 0 \pmod{p}$). Let $M_j, j = 1, 2, 3, \ldots, (n-1)$, be the $(p-1) * (p-n+j)$ matrix having as its first $(p-n)$ columns the columns of $M_0$ and as its last columns the first j columns of A. In general, $M_j (-U_{p-n+j}, -U_{p-n+j-1}, -U_{p-n+j-2}, \ldots, -U_1), j = 0, 1, 2, \ldots, (n-1)$, is congruent mod $p$ to the $(j+1)th$ column of $A$.

Then since the last column of $M_0$ is linearly dependent $(mod\ p)$ on the other columns of $M_0$, and the first column of A is linearly dependent $(mod\ p)$ on the columns of $M_0$, and the second column of A is linearly dependent $(mod\ p)$ on the columns of $M_1$, etc., there are at most $p - n - 1$ linearly independent columns of A. Since the rank of a matrix is the dimension of its column space, the rank of A is at most $p - 1 - n$. Then by Kronecker's version of König's theorem, the congruence $x^n + a_1 x^{n-1} + a_2 x^{n-2} + \ldots + a_n \equiv 0 \pmod{p}, 0 < x < p$, has at least $n$ roots. An $n$th degree mod p congruence has at most $n$ roots, therefore $x^n + a_1 x^{n-1} + a_2 x^{n-2} + \ldots + a_n \equiv 0 \pmod{p}, 0 < x < p$, has exactly $n$ roots.

## 2.5 Proof of Theorem (2.1) (the Part Using König's Theorem)

Now suppose $x^n + a_1 x^{n-1} + a_2 x^{n-2} + \ldots + a_n \equiv 0 \pmod{p}$, $0 < x < p$, has n roots. Then by Fermat's theorem, $A(x^{p-2}, x^{p-3}, x^{p-4}, \ldots, x^0) \equiv (0,0,0,\ldots,0) \pmod{p}$. Furthermore, by König's theorem, $p$ divides all $p-n$ rowed minors of $A$. One such minor is the determinant of the following matrix (denoted by $B_1$);

$$\begin{matrix} 0 & 0 & 0 & \ldots & 0 & 0 & \ldots & 0 & 0 & 1 & a_1 \\ 0 & 0 & 0 & \ldots & 0 & 0 & \ldots & 0 & 1 & a_1 & a_2 \\ 0 & 0 & 0 & \ldots & 0 & 0 & \ldots & 1 & a_1 & a_2 & a_3 \\ & & & & & \cdot & & & & & \\ 1 & a_1 & a_2 & \ldots & a_{n-1} & a_n & \ldots & 0 & 0 & 0 & 0 \\ a_1 & a_2 & a_3 & \ldots & a_n & 0 & \ldots & 0 & 0 & 0 & 0 \end{matrix}$$

(This is the matrix obtained from $A$ by deleting its first $n-1$ columns and last $n-1$ rows.) Therefore $B_1(y_{p-n-1}, y_{p-n-2}, y_{p-n-3}, \ldots, y_0) \equiv (0,0,0,\ldots,0) \pmod{p}$ where not all of $y_{p-n-1}, y_{p-n-2}, y_{p-n-3}, \ldots, y_0$ are congruent to zero mod $p$. Then $y_1 + a_1 y_0 \equiv 0 \pmod{p}$, $y_2 + a_1 y_1 + a_2 y_0 \equiv 0 \pmod{p}$, $y_3 + a_1 y_2 + a_2 y_1 + a_3 y_0 \equiv 0 \pmod{p}$, ..., $y_{p-n-1} + a_1 y_{p-n-2} + a_2 y_{p-n-3} + \ldots + a_n y_{p-2n-1} \equiv 0 \pmod{p}$ and hence $p$ does not divide $y_0$ (since otherwise $p$ would divide all of $y_{p-n-1}, y_{p-n-2}, y_{p-n-3}, \ldots, y_0$, a contradiction). Therefore $(y_1/y_0) + a_1 \equiv 0 \pmod{p}$. Also $U_1 + a_1 U_0 = 0$, therefore $(y_1/y_0) \equiv U_1 \pmod{p}$. Similarly $(y_2/y_0) + a_1 (y_1/y_0) + a_2 \equiv 0 \pmod{p}$ and $U_2 + a_1 U_1 + a_2 U_0 = 0$, therefore $(y_2/y_0) \equiv U_2 \pmod{p}$. By an induction argument,

$$(y_{p-n-1}/y_0) \equiv U_{p-n-1} \pmod{p}, (y_{p-n-2}/y_0) \equiv U_{p-n-2} \pmod{p},$$

$$(y_{p-n-3}/y_0) \equiv U_{p-n-3} \pmod{p}, \ldots, (y_1/y_0) \equiv U_1 \pmod{p}$$

and hence $B_1(U_{p-n-1}, U_{p-n-2}, U_{p-n-3}, \ldots, U_0) \equiv (0,0,0,\ldots,0) \pmod{p}$. The product of the last row of $B_1$ and $((U_{p-n-1}, U_{p-n-2}, U_{p-n-3}, \ldots, U_0)$ gives $a_1 U_{p-n-1} + a_2 U_{p-n-2} + \ldots + a_n U_{p-2n} \equiv 0 \pmod{p}$. Then since $U_{p-n} + a_1 U_{p-n-1} + a_2 U_{p-n-2} + \ldots + a_n U_{p-2n} = 0$, $U_{p-n} \equiv 0 \pmod{p}$.

Since $p$ divides all $p-n$ rowed minors of $A$, $p$ divides all $j$ rowed minors of $A$ where $j \geq p-n$. Therefore $p$ divides the determinant of the $(p-n+1)*(p-n+1)$ matrix

$$\begin{matrix} a_n & 0 & 0 & \ldots & 0 & 0 & \ldots & 0 & 0 & 1 & a_1 \\ 0 & 0 & 0 & \ldots & 0 & 0 & \ldots & 0 & 1 & a_1 & a_2 \\ 0 & 0 & 0 & \ldots & 0 & 0 & \ldots & 1 & a_1 & a_2 & a_3 \\ & & & & & \vdots & & & & & \\ 1 & a_1 & a_2 & \ldots & a_{n-1} & a_n & \ldots & 0 & 0 & 0 & 0 \end{matrix}$$

$$a_1 \; a_2 \; a_3 \; ... \; a_n \; 0 \; ... \; 0 \; 0 \; 0 \; 0$$

(This is the matrix obtained from $A$ by deleting its first $n-2$ columns and last $n-2$ rows.) $p$ divides the row 1, column 1 cofactor of this matrix (since $p$ divides all $p-n$ rowed minors of $A$), therefore $p$ divides the determinant of the matrix obtained from the above matrix by changing the row 1, column 1 element to a zero. Denote this matrix by $B_2$. Then $B_2(U_{p-n}, U_{p-n-1}, U_{p-n-2}, \ldots, U_0) \equiv (0, 0, 0, \ldots, 0)(mod\; p)$. The product of the last row of $B_2$ and $(U_{p-n}, U_{p-n-1}, U_{p-n-2}, \ldots, U_0)$ gives $a_1 U_{p-n} + a_2 U_{p-n-1} + a_3 U_{p-n-2} + \ldots + a_n U_{p-2n+1} \equiv 0(mod\; p)$. Then since $a_1 U_{p-n} + a_2 U_{p-n-1} + a_3 U_{p-n-2} + \ldots + a_n U_{p-2n+1} \equiv 0(mod\; p)$., $U_{p-n+1} \equiv 0(mod\; p)$.

The proofs that $U_{p-2} \equiv 0(mod\; p), U_{p-3} \equiv 0(mod\; p), U_{p-4} \equiv 0(mod\; p), \ldots, U_{p-n} \equiv 0(mod\; p)$ are similar. Finally, p divides the determinant of the $(p-1) \times (p-1)$ matrix

$$\begin{matrix} 0 & 0 & 0 & ... & 0 & & 0 & ... & 0 & 0 & 1 & a_1 \\ 0 & 0 & 0 & ... & 0 & & 0 & ... & 0 & 1 & a_1 & a_2 \\ 0 & 0 & 0 & ... & 0 & & 0 & ... & 1 & a_1 & a_2 & a_3 \\ & & & & & . & & & & & & \\ 1 & a_1 & a_2 & ... & a_{n-1} & a_n & ... & 0 & 0 & 0 & 0 \\ a_1 & a_2 & a_3 & ... & a_n & 0 & ... & 0 & 0 & 0 & 1 \end{matrix}$$

Denote this matrix by $B_n$. Then $B_n(U_{p-2}, U_{p-3}, U_{p-4}, \ldots, U_0) \equiv (0, 0, 0, \ldots, 0)(mod\; p)$. The product of the last row of $B_n$ and $(U_{p-2}, U_{p-3}, U_{p-4}, \ldots, U_0)$ gives $a_1 U_{p-2} + a_2 U_{p-3} + a_3 U_{p-4} + \ldots + a_n U_{p-n-1} + 1 \equiv 0(mod\; p)$. Then since $a_1 U_{p-2} + a_2 U_{p-3} + a_3 U_{p-4} + \ldots + a_n U_{p-n-1} = 0$, $U_{p-1} \equiv 1(mod\; p)$. Therefore $U_{p-1+i-n} \equiv U_{i-n} \; (mod\; p), i = 1, 2, 3, \ldots, n$.

**Theorem 2.6:** If $U_{p-1+i-n} \equiv U_{i-n} \; (mod\; p)$, $i = 1, 2, 3, \ldots, n$, then $U_{p-1+i-n} \equiv U_{i-n} \; (mod\; p), i = 1, 2, 3, \ldots$.

Suppose $U_{p-1+i-n} \equiv U_{i-n} \; (mod\; p), i = 1, 2, 3, \ldots, n$. $U_p + a_1 U_{p-1} + a_2 U_{p-2} + \ldots + a_n U_{p-n} = 0$, therefore $U_p + a_1 \equiv 0(mod\; p)$. Also $U_1 + a_1 U_0 = 0$, therefore $U_p \equiv U_1 \; (mod\; p)$, that is, $U_{p-1+i-n} \equiv U_{i-n} \; (mod\; p), i = n+1$. Similarly $U_{p+1} + a_1 U_p + a_2 U_{p-1} + \ldots + a_n U_{p-n+1} = 0$, therefore $U_{p+1} + a_1 U_1 + a_2 U_0 \equiv 0(mod\; p)$. Also $U_2 + a_1 U_1 + a_2 U_0 \equiv 0$, therefore $U_{p+1} \equiv U_2 \; (mod\; p)$, that is, $U_{p-1+i-n} \equiv U_{i-n} \; (mod\; p), i = n+2$. An induction argument gives $U_{p-1+i-n} \equiv U_{i-n} \; (mod\; p), i = 1, 2, 3, \ldots$. Therefore if $x^n + a_1 x^{n-1} + a_2 x^{n-2} + \ldots + a_n \equiv 0(mod\; p), 0 < x < p$, has $n$ roots, then $U_{p-1+i-n} \equiv U_{i-n} \; (mod\; p), i = 1, 2, 3, \ldots$.So $x^n + a_1 x^{n-1} + a_2 x^{n-2} + \ldots + a_n \equiv 0(mod\; p), 0 < x < p$, has $n$ roots if and only if $U_{p-1+i-n} \equiv U_{i-n} \; (mod\; p), i = 1, 2, 3, \ldots$.

## 2.6 Proofs of Theorems (2.2) and (2.4)

Denote $\sum \zeta_1^{e_1} \zeta_2^{e_2} \ldots \zeta_n^{e_n}$ where $\zeta_1, \zeta_2, \ldots \zeta_n$ are the roots of $\zeta^n + a_1 \zeta^{n-1} + a_2 \zeta^{n-2} + \ldots + a_n \equiv 0$ and the summation is over all combinations of non-negative integers $e_1, e_2, \ldots e_n$ such that $e_1 + e_2 + \ldots + e_n = i$ by $U'_i$. Then $U'_0, U'_1, U'_2, \ldots$ are defined and $U'_0 = 1$. If $i < 0$, let $U'_i = 0$. Then

**Theorem 2.7:** If $i \neq 0$, $U'_i = V_{i,1} + V_{i,2} + V_{i,3} + \ldots + V_{i,n}$.

$\zeta^n + a_1\zeta^{n-1} + a_2\zeta^{n-2} +\ldots+a_n = (\zeta - \zeta_1)(\zeta - \zeta_2)(\zeta - \zeta_3) \cdots (\zeta - \zeta_n)$, therefore $a_j, j = 1, 2, 3, \ldots, n$, equals $(-1)^j$ times the sum of all combinations of products of $\zeta_1, \zeta_2, \ldots \zeta_n$ taken j at a time (that is, $a_j = (-1)^j V_{j,j}$). If $i \geq 0$, each term in the summation giving $V_{i+j,j+k}$, $0 \leq k \leq i, k \leq n-j$, can be factored in $c(j+k,j)$ ways so that one factor is a term in the summation giving $U_i'$ and the other factor is a term in the summation giving aj. Conversely, if $i \geq 0$, every term in the summation giving $a_j U_i'$ is in one of $V_{i+j,j}, V_{i+j,j+1}, V_{i+j,j+2}, \ldots, V_{i+j,d}$ where $d = min(i+j, n)$. Therefore if i≥0,

$$a_j U_i' = (1)^j [c(j,j)V_{i+j,j} + c(j+1,j)V_{i+j,j+1} + c(j+2,j)V_{i+j,j+2}+\ldots+c(d,j)V_{i+j,d}].$$

Then $a_j U_i' = (1)^j [c(j,j)V_{i+j,j} + c(j+1,j)V_{i+j,j+1} + c(j+2,j)V_{i+j,j+2}+\ldots+c(n,j)V_{i+j,n}], j = 1, 2, 3, \ldots, n$. Denote the $n \times n$ matrix

$$\begin{matrix} -c(1,1) & -c(2,1) & -c(3,1) & \ldots & -c(n,1) \\ 0 & c(2,2) & c(3,2) & \ldots & c(n,2) \\ 0 & 0 & -c(3,3) & \ldots & -c(n,3) \\ & & & \vdots & \\ 0 & 0 & 0 & \ldots & (-1)^n c(n,n) \end{matrix}$$

by $T$. Therefore $T(V_{i,1}, V_{i,2}, V_{i,3}, \ldots, V_{i,n}) = (a_1 U_{i-1}', a_2 U_{i-2}', a_3 U_{i-3}', \ldots, a_n U_{i-n}')$. The sums of the columns of $T$ equal $-1$, therefore $-(V_{i,1} + V_{i,2} + V_{i,3}+\ldots+V_{i,n}) = a_1 U_{i-1}' + a_2 U_{i-2}' + a_3 U_{i-3}' + \ldots + a_n U_{i-n}'$ and hence $U_i' + a_1 U_{i-1}' + a_2 U_{i-2}' + a_3 U_{i-3}' + \ldots + a_n U_{i-n}' = 0, i = 1, 2, 3, \ldots$. Therefore $U_i' = U_i$ and Theorems (2.2) and (2.4) follow.

## 2.7 Proof of Theorem (2.5)

Element $(a,b), a \leq b$, of $T$ is $(-1)^a c(b,a)$ and element $(a,b)$ $a > b$, is zero, therefore element $(a,b), a \leq b$, of $T^2$ is $\sum(-1)^{k+a} c(k,a)c(b,k)$ where the summation is from $k = a$ to $k = b$, and element $(a,b), a > b$, is zero. If $a \leq k \leq b$, then $c(k,a)c(b,k) = c(b,a)c(b-a, k-a)$. Also, $\sum(-1)^{k+a} c(b,a)c(b-a, k-a)$ where the summation is from $k = a$ to $k = b$ equals $c(b,a)\sum(-1)^k c(b-a,k)$ where the summation is from $k = 0$ to $k = b-a$, and these summations equal 0 if $a < b$ or 1 if $a = b$.

Therefore $T^2 = I$ where $I$ is the $n \times n$ identity matrix. Then since $T(V_{i,1}, V_{i,2}, V_{i,3}, \ldots, V_{i,n}) = (a_1 U_{i-1}, a_2 U_{i-2}, a_3 U_{i-3}, \ldots, a_n U_{i-n})$, $T(a_1 U_{i-1}, a_2 U_{i-2}, a_3 U_{i-3}, \ldots, a_n U_{i-n}) = (V_{i,1}, V_{i,2}, V_{i,3}, \ldots, V_{i,n})$. Therefore $V_{i,j} = (-1)^j [c(j,j)a_j U_{i-j} + c(j+1,j)a_{j+1}U_{i-j-1} + c(j+2,j)a_{j+2}U_{i-j-2}+\ldots+c(n,j)a_n U_{i-n}], j = 1, 2, 3, \ldots, n$. (Note that this is the formula relating Fibonacci and Lucas numbers if $n = 2$, $a_1 = a_2 = -1$, and $j = 1$.)

## 2.8 Proofs of Remaining Theorems

Some previous assertions will now be proved.

$$\zeta_1^p + \zeta_2^p + \zeta_3^p+\ldots+\zeta_n^p \equiv (\zeta_1 + \zeta_2 + \zeta_3+\ldots+\zeta_n)^p \pmod{p}$$

(by properties of symmetric functions and the binomial coefficients) and

$$(\zeta_1 + \zeta_2 + \zeta_3 + \ldots + \zeta_n)^p \equiv (\zeta_1 + \zeta_2 + \zeta_3 + \ldots + \zeta_n)(mod\ p)$$

(by Fermat's theorem), therefore $\zeta_1^p + \zeta_2^p + \zeta_3^p + \ldots + \zeta_n^p \equiv (\zeta_1 + \zeta_2 + \zeta_3 + \ldots + \zeta_n)(mod\ p),$ that is;

**Theorem 2.8:** $V_{p,1} \equiv V_{1,1}\ (mod\ p)$.

Also, $V_{1,1} = -a_1$ and $-V_{p,1} = a_1 U_{p-1} + 2a_2 U_{p-2} + 3a_3 U_{p-3} + \ldots + na_n U_{p-n}$ (by Theorem (2.5)), therefore;

**Theorem 2.9:** $a_1(U_{p-1} - 1) + 2a_2 U_{p-2} + 3a_3 U_{p-3} + \ldots + na_n U_{p-n} \equiv 0(mod\ p)$.

Finally, Theorem (2.3) will be proved. If $U_{p-1+i-n} \equiv U_{i-n}\ (mod\ p), i = 1,2,3,\ldots,n$, then by the matrix equation obtained in the proof of Theorems (2.2) and (2.4), $V_{p,2} \equiv V_{p,3} \equiv V_{p,4} \equiv \ldots \equiv V_{p,n} \equiv 0(mod\ )$, that is $V_{p,1} \equiv V_{1,1}\ (mod\ p)\ i = 2,3,4,\ldots,n$. Conversely, if $V_{p,1} \equiv V_{1,1}\ (mod\ p)\ i = 2,3,4,\ldots,n$ and p does not divide $a_1 a_2 a_3 \ldots a_n$, then $U_{p-1} \equiv 1, U_{p-2} \equiv U_{p-3} \equiv U_{p-4} \equiv \ldots \equiv U_{p-n} \equiv 0(mod\ p)$, that is, $U_{p-1+i-n} \equiv U_{i-n}\ (mod\ p)$, $i = 1,2,3,\ldots,n$. Theorem (2.3) then follows from Theorem (2.1).

# 3. Conclusion

In this article we have find a new connection between $n^{th}$ degree polynomial mod $p$ congruence with $n$ roots and higher order Fibonacci and Lucas sequences. We have also discussed about the recent work that has been done in sequences and their connection to polynomial congruence and then find out new relations between particular recurrence relation and the congruence of the sequences.